\documentclass[12pt]{amsart}
\makeindex
\usepackage[latin1]{inputenc} \usepackage[T1]{fontenc} \usepackage{lmodern}
\usepackage{epsfig,graphicx,color,amsmath,a4wide,amssymb,amsfonts,amsbsy,xy,latexsym,epsf,pstricks,verbatim,times,hyperref}

\usepackage{pgf,color}
\usepackage{pgfarrows}

\xyoption{all}

\definecolor{Grey}{rgb}{.5,.5,.5}

\definecolor{Blue}{rgb}{.0,.0,0.9}
\definecolor{LightBlue1}{rgb}{.2,.4,0.9}
\definecolor{LightBlue2}{rgb}{.3,.5,0.9}
\definecolor{LightBlue3}{rgb}{.4,.6,0.9}
\definecolor{LightBlue4}{rgb}{.5,.7,.9}
\definecolor{LightBlue5}{rgb}{.6,.8,.9}
\definecolor{LightBlue6}{rgb}{.7,.9,.9}

\definecolor{Red}{rgb}{.9,.0,.0}
\definecolor{LightRed1}{rgb}{0.9,.2,.4}
\definecolor{LightRed2}{rgb}{0.9,.3,.5}
\definecolor{LightRed3}{rgb}{0.9,.4,.6}
\definecolor{LightRed4}{rgb}{.9,.5,.7}
\definecolor{LightRed5}{rgb}{.9,.6,.8}
\definecolor{LightRed6}{rgb}{.9,.7,.9}

\xyoption{all}
\usepackage[english]{babel}

\newcounter{noalgo}[section]
\newdimen\indentalgo
\newdimen\indentalgodec\indentalgo=0.0mm\indentalgodec=10mm

\newcommand{\If}{\advance\indentalgo by \indentalgodec {\bf if }}

\newcommand{\For}{\global\advance\indentalgo by \indentalgodec {\bf for }}

\newcommand{\Endindent}{\global\advance\indentalgo by -\indentalgodec}

\newdimen\decalage \decalage=0.5cm
\newcounter{algo} 

\def\<<{\leavevmode
  \raise0.28ex\hbox{$\scriptscriptstyle\langle\!\langle$}\nobreak
  \hskip -.6pt plus.3pt minus.2pt\,}
\def\>>{\,\nobreak\hskip -.6pt plus.3pt minus.2pt
  \raise0.28ex\hbox{$\scriptscriptstyle\rangle\!\rangle$}}

\def\<<{\leavevmode
  \raise0.28ex\hbox{$\scriptscriptstyle\langle\!\langle$}\nobreak
  \hskip -.6pt plus.3pt minus.2pt\,}
\def\>>{\,\nobreak\hskip -.6pt plus.3pt minus.2pt
  \raise0.28ex\hbox{$\scriptscriptstyle\rangle\!\rangle$}}

\newtheorem{definition}{D{e}finition}

\newtheorem{theorem}{Theorem}[section]

\newtheorem{lemma}{Lemma}[section]

\def\bK{{\mathbf L}}

\def\bF{{\mathbf F}}

\def\bK{{\mathbf K}}

\def\cP{{\mathbb P}}
\def\bA{{\mathbb A}}
\def\ZZ{{\mathbb Z}}

\def\bK{{\mathbf K}}
\def\bL{{\mathbf L}}

\def\Ker{{\mathop{\rm Ker}\nolimits}}

\def\Gal{{\mathop{\rm Gal}\nolimits}}

\def\div{\mathop{\rm{div}}\nolimits }
\def\Div{\mathop{\rm{Div}}\nolimits }

\def\ord{\mathop{\rm{ord}}\nolimits }

\def\Hom{\mathop{\rm{Hom}}\nolimits }
\def\Norm{\mathop{\rm{Norm}}\nolimits }

\newenvironment{myproof}[1][\myproofname]{\par
  \normalfont \topsep6pt\relax
  \trivlist
\item[\hskip\labelsep
  \itshape
  #1.]\ignorespaces
}{%
  \endtrivlist\hfill$\square$
}
\providecommand{\myproofname}{Proof}

\begin{document}

\title{Normal Bases using
1-dimensional Algebraic Groups}

\thanks{Research supported by the Simons Foundation via the PREMA project, and the Inria International Lab LIRIMA via the Associate team FAST.}

\author{Tony Ezome}
\address{Tony Ezome, Universit{\'e} des Sciences et Techniques de Masuku,
Facult{\'e} des Sciences, D{\'e}partement de math{\'e}matiques et informatique,
BP 943 Franceville, Gabon.}
\address{Tony Ezome,INRIA, LIRIMA, F-33400 TALENCE, FRANCE.}
\email{tony.ezome@gmail.com}

\author{Mohamadou Sall}
\address{Mohamadou Sall, Laboratoire d'Alg\`ebre, de Cryptologie, de G\'e om\'etrie Alg\'ebrique et Applications Universit\'e Cheikh Anta Diop de Dakar, BP 5005 Dakar Fann, S\'en\'egal}
\address{Mohamadou Sall, INRIA, LIRIMA, F-33400 TALENCE, FRANCE.}
\email{msallt12@gmail.com}

\maketitle

\maketitle

\begin{abstract}
This paper surveys and illustrates geometric
methods for constructing normal bases allowing
efficient finite field arithmetic. These bases are constructed
using the additive group, the multiplicative group and the Lucas torus.
We describe algorithms with quasi-linear complexity
to multiply two elements given in each one of the bases.
\end{abstract}


.

\vspace{.75cm}

\section{Introduction}

Consider two fields $\bK$ and $\bL$, such that $\bL$ is a 
degree $n$ cyclic extension of $\bK$. Denote 
by $\sigma$ a generator of the Galois group $\Gal(\bL/\bK)$.
A normal basis of
$\bL$ over $\bK$ is a basis
$(\theta,\sigma(\theta),\ldots,\sigma^{n-1}(\theta))$
generated by some $\theta$ in $\bL\!^*$. Such a
$\theta$ is called a normal
element of $\bL$ over $\bK$.
The normal basis theorem ensures that $\bL$ possesses
at least one normal element over $\bK$.

Let $\Theta=(\theta_i)_{0\le i\le n-1}$ be
an arbitrary basis of $\bL/\bK$. 
Given 
$$a=\sum_{i=0}^{n-1} a_i\theta_i\ \text{ and } \ 
b=\sum_{j=0}^{n-1} b_j\theta_j$$
in $\bL$, the sum $a+b$ is component-wise
and easy to implement.
The product $a\times b$  may be more difficult.
Let $\Gamma$ be a straight-line program computing the coordinates
of $a\times b$ in $\Theta$, from the coordinates of $a$ and $b$.
We assume that $\Gamma$ consists of additions,
subtractions, multiplications of a register by
a constant, and additions, subtractions, multiplications
between two registers. The complexity of
$\Gamma$ is the total number of such operations.
We define the complexity of $\Theta$ to be the minimal possible complexity
of a straight-line program computing the coordinates
of $a \times b$ from the ones of $a$ and $b$.
Let $t_{i,j}^{k}$ be coefficients in $\bK$ such that
\begin{equation}\label{eq:29} 
\theta_i\theta_j=\sum_{k=0}^{n-1} t_{i,j}^{k} \theta_k.
\end{equation}
Then
$$a\times b=\sum_{k=0}^{n-1} c_k(a,b)\theta_k,$$
where $c_k$ is a bilinear form on $\bL \times \bL$
defined by $$c_k(a,b) = \sum_{i,j}t_{i,j}^{k}a_ib_j.$$
Assume that $\Theta$ is a normal basis. Then every $x$
with vector coordinate $(x_0,x_1,\ldots,x_{n-1})$ in $\Theta$ is such that
$\sigma^{k}(x)$ has coordinate vector $(x_{-k},x_{-k+1},\ldots,x_{-k-1})$.
Since the coordinate vector of the product $\sigma^{n-k}(a)\times \sigma^{n-k}(b)$ is
equal to
$$(c_0(\sigma^{n-k}(a), \sigma^{n-k}(b)), c_1(\sigma^{n-k}(a), \sigma^{n-k}(b)),
\ldots, c_{n-1}(\sigma^{n-k}(a), \sigma^{n-k}(b))),$$
we have
$$ c_k(a,b)=c_0(\sigma^{n-k}(a),\sigma^{n-k}(b)).$$
This means that $c_k$ is obtained from $c_0$ by
a $k$-fold cyclic shift of coordinates of the variables involved.
Hence we define the \textit{weight}, denoted $w$, of a normal basis
to be the number of non-zero terms in the
form $c_{0}$. 
There is a straightforward
algorithm with complexity $2nw+n(w-1)$
for computing the coordinates of $a \times b$ from 
the ones of $a$ and $b$ in a normal
basis with weight $w$. The weight is sometimes called the complexity of the
normal basis, but we prefer to use a different terminology.
Using action of $\sigma$ on equation $( \ref{eq:29} )$, it is easily checked that
the weight of a normal basis $\Theta=(\theta_i)_{0\le i\le n-1}$ is
 also equal to the number of non-zero coefficients
 in the linear combinations
\begin{equation}\label{eq:30} 
\theta_0\theta_i=\sum_{j=0}^{n-1} t_{0,i}^{j} \theta_j ,  \text{ for }  \
1\le i \le n-1. 
\end{equation}
Mullin, Onyszchuk, Vanstone and Wilson \cite{Mullin-al} 
showed that the weight of any normal basis of $\bF\!_{q^{n}}$ over $\bF\!_q$
is greater than $2n-1$. This lower
bound is reached by the so-called optimal normal bases.
It is appropriate here to define
\textit{Gauss periods}.

\begin{definition} Let $q$ be a prime power. Let
 $n$ and $k$ be two integers such that
$r=nk+1$ is a prime number not dividing $q$. Denote by
$\mathcal{K}$ the unique subgroup of $(\ZZ/r\ZZ)^*$ of order $k$.
A Gauss period of type $(n,k)$ over $\bF\!_q$ is a sum of the form
$$\vartheta=\sum_{a \in \mathcal{K}} \theta^a$$
where $\theta$ is an arbitrary primitive $r$-th root of unity in
$\bF\!_{q^{nk}}$.
\end{definition}

It is easy to see that a Gauss period of type $(n,k)$ over $\bF\!_q$
lies in $ \bF\!_{q^n}$. Moreover,  a Gauss period 
of type $(n,k)$ over $\bF\!_q$ generates a normal basis of
$\bF\!_{q^{n}}/\bF\!_q$ if and only if $\gcd(e,n)=1$, where $e$
denotes the index of $q$ modulo $r=nk+1$ [see \cite{Gao-Gathen-Panario-Shoup},
or \cite{Wassermann1}].
Optimal normal bases occur when $k = 1$ or $k = 2$.
In fact normal bases with low weight and low complexity are usually constructed using
\textit{Gauss periods} [see \cite{Ash-Blake-Vanstone}, 
\cite{Christopoulou-Garefalakis-Panario-Tomson1},
\cite{Gao-Lenstra}, \cite{Gao-Gathen-Panario-Shoup},
\cite{Liao-You}, \cite{Wan-Zhou}]. 
But Gao described in [\cite{Gao-Phd}, chapter 5] another way to construct
normal bases with low weight. The Lucas torus
and its isogenies play an important, though implicit, role in Gao's construction.
Further, Couveignes and Lercier constructed normal bases using elliptic curves \cite{Couveignes-Lercier2}.
The resulting elliptic
normal bases allow quasi-linear finite field
arithmetic.
Our work is concerned by efficient normal bases constructed with
1-dimensional algebraic groups. It is known that there are only
three 1-dimensional connected affine algebraic
groups over a perfect field (up to isomorphism): the additive group
$\mathbf{G}_a$, the multiplicative $\mathbf{G}_m$ and the Lucas torus
$\mathbb{T}_\alpha$ [see \cite{milne2006}, chapter II, section 3, page 54].
Moreover, any connected projective algebraic
group of dimension 1 is an elliptic curve.
So a 1-dimensional connected algebraic group $G$ over a perfect field
is either $\mathbf{G}_a$, $\mathbf{G}_m$, $\mathbb{T}_\alpha$
or an elliptic curve. In this paper, we consider normal bases constructed
from $\mathbf{G}_a$, $\mathbf{G}_m$, and $\mathbb{T}_\alpha$.
We adapt the construction of elliptic normal bases proposed by Couveignes and Lercier
to these contexts. That results in natural and efficient algorithms.
We prove the following theorems:

\begin{theorem}\label{theorem:1}
Let $\bK$ be a field with characteristic $p>0$.
Let $a\in \bK$ be an element which
does not lie in $\{x^p-x|x\in \bK\}$.
Then $\bL:=\bK[X]/(X^p-X-a)$ is a degree $p$ cyclic extension of
$\bK$, and there is a normal basis of $\bL$ over 
$\bK$ with weight $\le 3p-2$ and complexity
$O(p(\log p)(\log|\log p|))$.
\end{theorem}

\begin{theorem}\label{theorem:Kummer}
Let $\bK$ be a field with characteristic $p>0$.
Assume that $m\ge 2$ and $n\ge 2$ are two integers
such that $mn$ is
prime to $p$. Assume that $\bK$ possesses
a primitive $mn$-th root of unity.
Let $a\in \bK$ be a non-zero element
such that the order of the class of $a$
in $\bK^*/\bK^{*n}$ is equal to $n$.
Then $\bL:=\bK[X]/(X^n-a)$ is a degree $n$ cyclic extension of
$\bK$, and there is a normal basis of $\bL$ over 
$\bK$ with weight $\le 3n-2$ and complexity
$O(n(\log n)(\log|\log n|))$.
\end{theorem}

If $\bK=\bF\!_q$ is a finite field, then the first
two requirements in theorem \ref{theorem:Kummer}
are equivalent to saying that
$mn$ divides $q-1$. A sufficient condition for the last
requirement in theorem \ref{theorem:Kummer}, in case $\bK=\bF\!_q$,
is to take for $a$ a generator of $\bF{\!_q}\!^*$.
Gao constructed low weight normal bases of $\bF\!_{q^n}/\bF\!_{q}$
in [\cite{Gao-Phd}, chapter 5, section 5.3] by using irreducible
polynomials of degree $n$ which  divide $cX^{q+1}+dX^q-aX-b \in \bF\!_q[X]$,
where $c\ne 0$, $ad-bc\ne 0$, and $n$ divides $q-1$.
The weights of the resulting normal bases
have the same upper bound as the one in theorem \ref{theorem:Kummer},
and the one in theorem \ref{theorem:1} when $n$ is equal to the characteristic of 
$\bF\!_q$.

\begin{theorem}\label{theorem:2}
Let $q$ be a prime power and $n$
a non-trivial divisor of $q+1$. 
Then there exists a normal basis $\Theta$ of $\bF\!_{q^n}$
over $\bF\!_q$ with complexity $O(n(\log n)(\log|\log n|))$.
\end{theorem}

Note that elliptic normal bases of $\bF\!_{q^n}/\bF\!_{q}$ 
constructed by Couveignes and Lercier \cite{Couveignes-Lercier2} have
complexity $O(n(\log n)^2(\log|\log n|))$ when they exist.
Gao, von zur Gathen, Panario and Shoup showed
\cite{Gao-Gathen-Panario-Shoup} that fast multiplication methods
(like FFT) can be adapted to normal bases 
of $\bF\!_{q^n}/\bF\!_q$ constructed with Gauss periods.
They proved that the  complexity of a normal basis of $\bF\!_{q^n}/\bF\!_q$
generated by a Gauss period of type $(n,k)$ is
equal to $O(nk(\log nk)(\log|\log nk|))$.

\subsection*{Plan}
In section 2 we prove theorem 1 using the additive group.
In section 3 we use the multiplicative group to prove theorem 2.
In sections 4 we explain how the Lucas Torus can be used to prove
theorem 3. At the end of each section, we give detailed examples.

\subsection*{Acknowledgments}
We thank Jean-Marc Couveignes
for his comments on early versions of this work.
The first author acknowledges 
the International Centre for Theoretical Physics (ICTP) and
the Mathematisches Forschungsinstitut Oberwolfach (MFO)
for their hospitality.

\section{Constructing normal bases with the additive group}

Consider a field $\bK$ with characteristic $p>0$.
We denote by $\overline{\bK}$ an algebraic closure
of $\bK$. We identify the additive group $\mathbf{G}_a$
over $\bK$ with the affine line
$\bA^1$ over $\bK$ endowed with the $x$-coordinate.
Any point $P$ in $\mathbf{G}_a$
is given by its $x$-coordinate.
The unit element $O_{\mathbf{G}_a}$ has $x$-coordinate equal to $0$.
The group law $\oplus_{\mathbf{G}_a}$ is defined by
$$x(P_1\oplus_{\mathbf{G}_a}P_2)=x(P_1)+x(P_2).$$

\subsection{Specializing isogenies of the additive group}\label{section:additive-group}

The $\bF\!_p$-rational points of  $\mathbf{G}_a$ form 
a cyclic subgroup of $\mathbf{G}_a(\bK)$.
Let $I : \mathbf{G}_a 
\rightarrow \mathbf{G}_a$ be the quotient isogeny
of $\mathbf{G}_a$ by $\mathbf{G}_a(\bF\!_p)$.
In terms of $x$-coordinates, $I$ is given by
$$x(I(P))=x(P)^p-x(P).$$
Let $a$ be a $\bK$-rational point in $\mathbf{G}_a$
outside of the image $I(\mathbf{G}_a(\bK))$.
Then the subfield $\bL=\bK(I^{-1}(a))$ of $\overline{\bK}$
is a cyclic extension of $\bK$ with degree $p$.
Indeed, fix $b$ in $I^{-1}(a)$
and denote by $\theta$ the $x$-coordinate of $b$.
Since $\bL$ is the splitting field of the separable polynomial
$$X^p-X-x(a)=(X-\theta)(X-(\theta-1)\cdots (X-(\theta-(p-1)))
 \in \bK[X],$$
it is a normal and separable extension of $\bK$.
If $P(X) \in \bK[X]$ is an irreducible factor of $X^p-X-x(a)$
with degree $1\le r<p$, then  $r\theta \in \bK$. 
There exist $u,v \in \ZZ$ such that $ur+vp=1$. So $\theta= ur\theta$
lies in $\bK$, contradicting the assumption that $a=I(b)$
does not lie in $I(\mathbf{G}_a(\bK))$.
Hence $X^p-X-x(a)$ is irreducible over $\bK$, and $\bL$ is a degree
$p$ Galois extension of $\bK$.              
The Galois group $\Gal(\bL/\bK)$ is made of $\bK$-automorphism
$\mathfrak{a}_k$ such that $\mathfrak{a}_k(\theta)=\theta+k$, for
$0\le k\le p-1$. So $\Gal(\bL/\bK)$ is generated by $\mathfrak{a}_1$.

We set $t= \mathfrak{a}_1(b)\ominus_{\mathbf{G}_a}b$.
Then the fiber of $I$ above $a$ is given by
\begin{equation}\label{eq:1}
I^{-1}(a)=[b]+[b\oplus_{\mathbf{G}_a}t]+\cdots+[b
\oplus_{\mathbf{G}_a}(p-1)t)].
\end{equation}

The additive group is an open subset
of the projective line $\cP^1$. Consider the divisor
$$D=[O_{\mathbf{G}_a}]+[t]+[2t]+
\cdots+[(p-1)t]-[\infty] \in \Div ( \cP^1).$$
The linear space
\begin{equation}\label{eq:2}
\mathcal{L}=\mathbf{H}^0(\cP^1,\mathcal{O}_{\cP^1}(D))
\end{equation}
has dimension $p$ over $\overline{\bK}$.
The translation $\tau : P\mapsto P\ominus_{\mathbf{G}_a} t$
 is an automorphism of
$\mathbf{G}_a$ which extends to the whole $\cP^1$
 by setting $\tau(\infty)=\infty$.
The divisor $D$ is invariant by $\tau$,
then so is the space $\mathcal{L}$.
For $0\le k\le p-1$, the functions 
$$\frac{1}{x\circ \tau^{k}} \in 
\bK(\mathbf{G}_a)$$
lie in $\mathcal{L}$.
Examination of
poles shows that they are linearly independent. So the system
 $$(
\frac{1}{x},\frac{1}{x-1},\ldots,\frac{1}{x-(p-1)})$$ is a basis of
$\mathcal{L}$ invariant by $\tau$.  
Evaluation at $b$ results
in a normal basis
\begin{equation}\label{eq:3}
\Theta=(\frac{1}{\theta},\frac{1}{\theta-1},\ldots,
\frac{1}{\theta-(p-1)})
\end{equation}
 of $\bL$ over $\bK$. 
Indeed, let $\lambda_0,\ldots,\lambda_{p-1}$ be
scalars in $\bK$ such that
$$\sum_{k\in \bF\!_p}\lambda_k\frac{1}{\theta-k}=0.$$
Then the function $\sum_{k\in \bF\!_p}\lambda_k\frac{1}{x-k}
=\frac{F(x)}{x(x-1) \cdots(x-(p-1))}$ cancels at $b$, where
$F(X)$ is a polynomial in $\bK[X]$.
If $F(X)$ is non-zero, then its degree
is at most $p-1$. But $F$ cancels at $\theta$
and also at all its $p$ conjugates,
this is impossible. So $F$ is the zero polynomial. 
All $\lambda_k$ are $0$ because
the functions $\frac{1}{x-k}$ are linearly independant.
Recall that $\Gal(\bL/\bK)=\{\mathfrak{a}_k|
0\le k\le p-1\}$ is generated by $\mathfrak{a}_1$
which satisfies $\mathfrak{a}_1(\theta)=\theta+1$.
So $$\mathfrak{a}_1(\frac{1}{\theta-k})=
\frac{1}{\theta-(k-1)},  \text{ for  } k \text{ in }  \bF\!_p.$$
We conclude that $\Theta$ is a normal basis.
If $i\ne 0$, then
$$\frac{1}{\theta}\times \frac{1}{\theta-i}=\frac{1}{i}(\frac{1}{\theta-i}
-\frac{1}{\theta}).$$
Since the weight can be defined
using linear combinations in equation $(\ref{eq:30})$,
we conclude that the weight of $\Theta$
is at most $3p-2$.

\subsection{Complexity}\label{section:complexity-additive-group}

In this section we describe an FFT-like algorithm which computes
the product of two elements of $ \bL$
in the normal basis $\Theta$ defined in  $(\ref{eq:3})$.
We adapt the construction proposed by Couveignes and Lercier
in [\cite{Couveignes-Lercier2}, section 4.3] to our context.\\

\noindent \textbf{Notation:} 
Let $\overrightarrow{\alpha}=(\alpha_k)_{0\le k\le p-1}$ and 
$\overrightarrow{\beta}=(\beta_k)_{0\le k\le p-1}$ be two vectors in
$\bK^{p}$. We denote by
$\overrightarrow{\alpha}\star_k \overrightarrow{\beta}
=\sum_{0\le i\le p-1} \alpha_i \beta_{k-i}$
the $k$-th component of the convolution product.
We denote by $\sigma(\overrightarrow{\alpha})=
(\alpha_{k-1})_{k}$ the cyclic shift of 
$\overrightarrow{\alpha}$. We denote by 
$\overrightarrow{\alpha} \diamond \overrightarrow{\beta}
=(\alpha_k \beta_k)_{k} $
the component-wise product and by
$\overrightarrow{\alpha}\star \overrightarrow{\beta}
=(\overrightarrow{\alpha}\star_k
\overrightarrow{\beta})_k$ the convolution product.\\

\noindent \textbf{Reduction and evaluation.} 
We fix $u=\frac{1}{x}, \text{ and }
\theta_0= \frac{1}{\theta}.$ For $0 \le k \le p-1$, we set
$$ u_k=u \circ  \tau^k,  \theta_k=
\mathfrak{a}_{-k}(\theta_0), \xi_k=\theta_k^2.$$

We also set $\xi_0=\sum_{k\in \bF\!_p} \imath_k
\theta_k$ and $\overrightarrow{\imath}=(\imath_k)_{0\le k\le 
p-1}$. We want to reduce a linear combination of the $\xi$'s into
a linear combination of the $\theta$'s. We have
$$\xi_i=\mathfrak{a}_1^{-i}(\xi_0)=\sum_{k\in \bF\!_p} \imath_k
\theta_{k+i}=\sum_{k\in \bF\!_p} \imath_{k-i}
\theta_k  \text{ for }\le i\le 
p-1.$$
Let 
$\overrightarrow{\alpha}=(\alpha_i)_{0\le i\le 
p-1}$ and $\overrightarrow{
\beta}=(\beta_j)_{0\le j
\le p}$ be two vectors in $\bK^p$ such that
$$\sum_{i\in \bF\!_p} \alpha_i\xi_i=\sum_{j \in \bF\!_p}
\beta_j\theta_j.$$
Since
$$
\sum_{i} \alpha_i\xi_i  =\sum_{i}
\alpha_i\sum_{k} \imath_{k-i}\theta_k
 =\sum_{k}
\theta_k\sum_{i} \alpha_i\imath_{k-i}
=\sum_{k} (\overrightarrow{\imath}\star_k
\overrightarrow{\alpha})\theta_k.
$$
We have
\begin{equation}\label{eq:4}
\beta_j=\overrightarrow{\imath}\star_k
\overrightarrow{\alpha}, \text{ that is }
\overrightarrow{\beta}=\overrightarrow{\imath} \star
 \overrightarrow{\alpha}.
\end{equation}

We now focus on evaluation of some functions
in $\bK(\mathbf{G}_a)$. Let $R$ be a point in $\mathbf{G}_a(\bK)$
outside of the subgroup generated by $t$.
We want to evaluate $f=\sum_{i\in \ZZ/p\ZZ} \alpha_iu_i$ 
at $R+jt$ for $0\le j \le p-1$.
We have
$$
f((R+jt)=\sum_{i\in \bF\!_p}
\alpha_i u_i(R+jt)=\sum_{i\in \bF_p}
\alpha_i u_0(R+(j-i)t)=\overrightarrow{\alpha}\star_j
\overrightarrow{u}_{\!\!R}
$$
where $\overrightarrow{u}_{\!\!R}=(u_0(R+kt))_{k\in \bF\!_p}$.
So the evaluation of $f$ is given by the convolution product
\begin{equation}\label{eq:5}
\overrightarrow{u}_{\!\!R}\star \overrightarrow{\alpha}.
\end{equation}

Similarly, for $f=\sum_{i\in \ZZ/n\ZZ} \alpha_iu^2_i$
we have 
$$
f(R+jt)=\sum_{i\in \bF\!_p}
\alpha_i u^2_i(R+jt)=\sum_{i\in \bF\!_p}
\alpha_i u^2_0(R+(j-i)t)=\overrightarrow{\alpha}\star_j
\overrightarrow{w}_{\!\!R}
$$
where $\overrightarrow{w}_{\!\!R}=(u^2_0(R+kt))_{k\in \bF\!_p}$.
So the evaluation of $f$ is given by the convolution product
\begin{equation}\label{eq:6}
\overrightarrow{w}_{\!\!R}\star \overrightarrow{\alpha}.
\end{equation}

\noindent \textbf{Interpolation.}  
The evaluation map
$f  \mapsto (f(R+jt))_{j\in \bF\!_p}$ is a bijection from 
the linear space $\mathcal{L}$ onto $\bK^p$.
Indeed two functions $f_1,f_2$ in $\mathcal{L}$ have the same evaluation
if and only if the function $f_1-f_2$
cancels at $R,R+t, \ldots, R+(p-1)t$ and $\infty$. But $f_1-f_2$
 has at most $p$ poles. So $f_1=f_2$.

Given a vector $\overrightarrow{\beta}=(\beta_0,\ldots,\beta_{p-1})
$ in $\bK^p$, we can compute the function $f$ in $\mathcal{L}$
such that $f(R+jt)=\beta_j$ by inverting the evaluation
map. That corresponds to the inverse $\overrightarrow{u}_{\!\!R}^{-1}$
of $\overrightarrow{u}_{\!\!R}$ for the convolution
product.

\noindent \textbf{An efficient multiplication algorithm.}
We want to compute the coordinates in $\Theta$ of the product
$$\left( \sum_{i\in \bF\!_p} \alpha_i\theta_i \right)  \times
\left( \sum_{i\in \bF\!_p}\beta_j\theta_j\right).$$
 Define the functions
$$\begin{array}{rl}
 & A=\sum_{i} \alpha_iu_i, B=\sum_{i} \beta_iu_i,\\
 & \\
& C=\sum_{i} \alpha_i\beta_iu_i^2,\\
 & \\
& D=AB-C.
\end{array}
$$
The product we want to compute is $A(b)B(b)=C(b)+D(b)$.\\
If $ i,j\in \bF_p$ are such  that $i\ne j$, then
 $u_iu_j$ lies in $\mathcal{L}$.
So 
$$(\sum_i\alpha_iu_i)(\sum_i\alpha_iu_i)=\sum_i\alpha_i\beta_iu_i^2 \bmod 
\mathcal{L}, $$
that is $D$ is in $\mathcal{L}$.
From equation $(\ref{eq:4})$, we deduce that
the coordinates in $\Theta$ of $C(b)$ are given by the vector
$$\overrightarrow{\imath}\star (\overrightarrow{\alpha}\diamond
\overrightarrow{\beta}).$$

From equation $(\ref{eq:5})$, the evaluation of $A$ at the points
$(R+jt)_j$ is given by
$\overrightarrow{u}_{\!\!R}\star \overrightarrow{\alpha}$.
The evaluation of $D$ at theses points is
$$(\overrightarrow{u}_{\!\!R}\star \overrightarrow{\alpha})
\diamond (\overrightarrow{u}_{\!\!R}\star \overrightarrow{\beta})
- \overrightarrow{w}_{\!\!R}\star (\overrightarrow{\alpha}
\diamond \overrightarrow{\beta}).$$
We get the coordinates of $D$ in 
the basis $(u_0,\ldots,u_{p-1})$ by applying $\overrightarrow{u}_{\!\!R}^{-1}$
on the left to this vector. These are also the coordinates
of $D(b)$ in the basis $\Theta$.

Altogether, the coordinates in $\Theta$ of
the product $(\sum_{i\in \bF_p} \alpha_i\theta_i)
\times (\sum_{j \in \bF\!_p}\beta_j\theta_j)$ are given by
$$ \overrightarrow{\imath}\star (\overrightarrow{\alpha}\diamond
\overrightarrow{\beta}) +
\overrightarrow{u}_{\!\!R}^{-1}\star \left((\overrightarrow{u}_{\!\!R}\star \overrightarrow{\alpha})
\diamond (\overrightarrow{u}_{\!\!R}\star \overrightarrow{\beta})
- \overrightarrow{w}_{\!\!R}\star (\overrightarrow{\alpha}
\diamond \overrightarrow{\beta})\right).$$

That consists in 5 convolution
products, 2 component-wise products, 1 addition
and 1 subtraction between vectors in $\bK^p$.

Each convolution product can be computed at the expense
of $O(p\log p |\log|\log p||)$ operations in $\bK$ using algorithms
due to Sch\"{o}nhage and Strassen \cite{schonhage-Strassen},
Sch\"{o}nhage \cite{schonhage}, Cantor and Kaltofen \cite{Cantor-Kaltofen}
(see [\cite{Gathen-Gerhard}, section 8.3] for a survey).
So the above multiplication algorithm has complexity 
$O(p(\log p)|\log|\log p||)$.

\subsection{Example}

Take $p=5$ and $\bK= \bF_5[X]/(X^3 +3X +2)$.
We set $\epsilon =X\bmod X^3+3X+2$.
We denote by $R$ the point in $\mathbf{G}_a$ with 
coordinate $x(R)=\epsilon$.
The point $a \in  \mathbf{G}_a(\bK)$ with coordinate $x(a)=1$ does not lie in
the image $I(\mathbf{G}_a(\bK))$ because
$Y^5-Y-1$ is relatively prime to $Y^{125}-Y$.
So $Y^5-Y-1$ is irreducible over $\bK$. We set
$\bL=\bK[Y]/(Y^5-Y-1)$ and $\theta=Y\bmod Y^5-Y-1$.
So $\Theta=(\theta_k)_{0\le k \le 4}$ is a normal basis
of $\bL/\bK$.
The weight of $\Theta$ is equal to $13$, according
to  following equations :

$$\theta_0^2 = 4\theta_0+4\theta_1+2\theta_2
+3\theta_3+\theta_4, \quad \theta_0\times\theta_1= -\theta_0+\theta_1 $$

$$\theta_0\times\theta_2= \frac{1}{2}(-\theta_0+\theta_2),
\quad \theta_0\times\theta_3= \frac{1}{3}(-\theta_0+\theta_3),
\quad \theta_0\times\theta_4= \frac{1}{4}(-\theta_0+\theta_4).$$

Now we  compute the coordinates 
 in $\Theta$ of the product 
$$\left( \sum_{i\in \bF_5} \alpha_i\theta_i \right)  \times
\left( \sum_{i\in \bF_5}\beta_j\theta_j\right) \text{  where  }
\overrightarrow{\alpha}=(1,3,1,1,2) \text{ and } \overrightarrow{\beta}
=(2,1,1,4,2).$$

We know that $\overrightarrow{\imath} =(4, 4, 2, 3, 1)$.
We compute
$$\overrightarrow{u}_{\!\!R}=(2\epsilon^{2}+1,4\epsilon^{2}+4\epsilon+1,
4\epsilon^2+3\epsilon+3,3\epsilon^2+4\epsilon+1,3\epsilon^2+2\epsilon+2),$$
$$  \overrightarrow{u}_{\!\!R}^{-1}=
(3\epsilon +4,   2\epsilon^2 + \epsilon + 4,  2\epsilon^2 + 
4\epsilon + 2,  3\epsilon,  3\epsilon^2 + \epsilon + 4), $$
$$\overrightarrow{w}_{\!\!R}=
(2\epsilon^2 + 2\epsilon + 1,  \epsilon^2 + 2,  4\epsilon + 1, 
 3\epsilon + 3,  4\epsilon^2 + 4\epsilon).$$
So
$$
\begin{array}{l}
 \overrightarrow{\imath}\star (\overrightarrow{\alpha}\diamond
\overrightarrow{\beta})  =  (3, 1, 1,1, 0)   \\
  \overrightarrow{u}_{\!\!R}\star \overrightarrow{\alpha}   = 
 (\epsilon^2 + \epsilon + 3,  4\epsilon^2 + \epsilon + 3,
 2\epsilon^2 + 1,  2\epsilon^2 + \epsilon + 1,  4\epsilon^2 + \epsilon + 1)   \\
  \overrightarrow{u}_{\!\!R}\star \overrightarrow{\beta}  = 
 (4\epsilon^2 + \epsilon + 4, 3\epsilon^2 + 2\epsilon,
  2\epsilon^2 + \epsilon + 3, 3\epsilon^2 + 4\epsilon + 4,
 3\epsilon^2 + 2\epsilon + 4)  \\
  \overrightarrow{w}_{\!\!R}\star (\overrightarrow{\alpha}
\diamond \overrightarrow{\beta})  = (2, 2\epsilon^2 + 3\epsilon + 3,
 \epsilon^2 + 3\epsilon + 1, 2\epsilon, 4\epsilon + 2)\\
 (\overrightarrow{u}_{\!\!R}\star \overrightarrow{\alpha}) \diamond 
 (\overrightarrow{u}_{\!\!R}\star \overrightarrow{\beta})   =
 (4\epsilon + 2, 4\epsilon + 3, \epsilon^2 + 2\epsilon + 4,
  2\epsilon^2 + 3\epsilon + 2, 4\epsilon + 2) .
\end{array}$$
Therefore
$$ 
\begin{array}{rl}
 \overrightarrow{u}_{\!\!R}^{-1}\star \left( ( \overrightarrow{u}_{\!\!R}\star
 \overrightarrow{\alpha}) \diamond ( \overrightarrow{u}_{\!\!R}\star
 \overrightarrow{\beta}) - \overrightarrow{w}_{\!\!R}\star
 (\overrightarrow{\alpha} \diamond \overrightarrow{\beta})\right)
 & =  (2, 2, 4, 4, 3).
\end{array}$$

Finally we get
$$\left( \sum_{i\in \bF_5} \alpha_i\theta_i \right)  \times
\left( \sum_{i\in \bF_5}\beta_j\theta_j\right) = 3\theta_1+3\theta_4.$$


\section{The multiplicative group case}\label{section:multiplicative-group}

Let $\bK$ be a field with characteristic $p > 0$.
Consider the affine line $\bA^1$ over $\bK$ endowed with the $x$-coordinate.
We identify the multiplicative group $\mathbf{G}_m$ with the open
subset $\{ x\ne 0 \}$ of $\bA^1$. Any point $P$ in $\mathbf{G}_m$ is given by
its $x$-coordinate.
The unit element $O_{\mathbf{G}_m}$ has $x$-coordinate equal to $1$.
The group law $\oplus_{\mathbf{G}_m}$ is defined by
$$x(P_1\oplus_{\mathbf{G}_m}P_1)=x(P_1)\times x(P_2).$$

\subsection{Specializing isogenies of the multiplicative group}\label{section:multiplicative-group}

Let $m\ge 2$ and $m\le 2$ be two integers
such that $mn$ is prime to $p$.
We assume that $\bK$ contains a primitive $mn$-th root of unity
which we denote by $\zeta_{mn}$. If $\bK=\bF\!_q$ is a finite field, 
 this is equivalent to saying that
$mn$ divides $q-1$. We set $\zeta_{n}=(\zeta_{mn})^m$.
The $n$-torsion $\mathbf{G}_m[n]$ is a cyclic subgroup
of order $n$ of $\mathbf{G}_m(\bK)$.
Let $I : \mathbf{G}_m 
\rightarrow \mathbf{G}_m$ be quotient isogeny
of $\mathbf{G}_m$ by $\mathbf{G}_m[n]$.
This is the multiplication by $n$ isogeny.
Let $a$ be a $\bK$-rational point in $\mathbf{G}_m$
such that $a\bmod I(\mathbf{G}_m(\bK))$ has order
$n$ in $\mathbf{G}_m(\bK)/I(\mathbf{G}_m(\bK))$.
In case $\bK=\bF\!_q$ is a finite field, take for $a$ 
a generator of $\bF\!_{q}^*$ is a sufficient condition for this last
requirement. 
In any case, the subfield 
$\bL=\bK(I^{-1}(a))$ of $\overline{\bK}$
is a cyclic extension of $\bK$ with degree $n$.
Indeed, fix $b$ in $I^{-1}(a)$
and denote by $\theta$ the $x$-coordinate of $b$.
Since $\bL$ is the splitting field of the separable polynomial
$$X^n-x(a)=(X-\theta)(X-\zeta_n \theta)(X-\zeta^2_n \theta)\cdots
 (X-\zeta^{n-1}_n \theta) \in \bK[X], $$
it is a normal and separable extension of $\bK$.
Let $P(X) \in  \bK[X]$ be a monic  irreducible factor of $X^n-x(a)$
with degree $r\le n$ such that $r$ is the smallest element
among the degrees of irreducible polynomials dividing $X^n-x(a)$.
Then $r$ is the smallest integer $>0$ such that 
$\theta^{r}  \in  \bK$. This means that $r$ is the smallest 
interger $>0$ such that $[r](a)  \in I(\mathbf{G}_m(\bK))$.
Since $a\bmod I(\mathbf{G}_m(\bK))$ has order
$n$ in $\mathbf{G}_m(\bK)/I(\mathbf{G}_m(\bK))$, we have
$r=n$ and $P(X)= X^n-x(a)$.
Hence $\bL$ is a degree $n$ Galois extension of $\bK$.
The Galois group $\Gal(\bL/\bK)$ is made of $\bK$-automorphisms
$\mathfrak{a}_k$ such that $\mathfrak{a}_k(\theta)=\zeta^{k}_n \theta$,
for $0\le k \le  n-1$. So $\Gal(\bL/\bK)$ is generated by $\mathfrak{a}_1$.
We set $t=\mathfrak{a}_1(b)\ominus_{\mathbf{G}_m} b$.
Then the fiber of $I$ above $a$ is given by
\begin{equation}\label{eq:7}
I^{-1}(a)=[b]+[b\oplus_{\mathbf{G}_m}t]+\cdots+[b
\oplus_{\mathbf{G}_m}(n-1)t)].
\end{equation}

The multiplicative group is an open 
subset of the projective line $\cP^1$.
Consider the divisor
$$D=[O_{\mathbf{G}_m}]+[t]+[2t]+
\cdots+[(n-1)t]-[\infty]\in \Div ( \cP^1).$$ 

The linear space
\begin{equation}\label{eq:8}
\mathcal{L}=\mathbf{H}^0(\cP^1,\mathcal{O}_{\cP^1}(D))
\end{equation}
has dimension $n$.
The translation $\tau : P\mapsto P\ominus_{\mathbf{G}_m}t$ is an automorphism of
$\mathbf{G}_m$ which extends to the whole $\cP^1$
($\tau(\infty)=\infty$).
The divisor $D$ is invariant by $\tau$,
then so is the space $\mathcal{L}$.
The functions 
$$\frac{1}{x-1},\frac{1}{\zeta_n^{-1} x-1},\ldots,\frac{1}{\zeta_n^{-(n-1)}x-1} \in 
\bK(\mathbf{G}_m)$$
lie in $\mathcal{L}$.
Examination of
poles shows that they are linearly independent. So they form a basis of
$\mathcal{L}$ invariant by $\tau$. 
Evaluation of these functions at $b$,
results in a normal basis
\begin{equation}\label{eq:9}
\Theta=(\frac{1}{\theta-1},\frac{1}{\zeta^{-1}\theta-1},\ldots,
\frac{1}{\zeta^{-(n-1)}\theta-1})
\end{equation}
 of $\bL$ over $\bK$.
Indeed, let $\lambda_0,\ldots,\lambda_{n-1}$ be
scalars in $\bK$ such that
$$\sum_{0\le k \le n-1}\lambda_k\frac{1}{\zeta^{-k}\theta-1}=0.$$
Then the function $\sum_{0\le k \le n-1}\lambda_k\frac{1}{\zeta_n^{-k}x-1}
=\frac{F(x)}{(x-1)(\zeta_n^{-1} x-1) \cdots(\zeta_n^{-(n-1)}x-1)}$ cancels at $b$, where
$F(X)$ is a polynomial in $\bK[X]$.
If $F(X)$ is non-zero, then its degree is 
 at most $\le n-1$. But $F$ cancels at $\theta$
and also at all its $n$ conjugates,
this is impossible. So $F$ is the zero polynomial. 
All $\lambda_k$ are $0$ because
the functions $\frac{1}{\zeta_n^{-k}x-1}$ are linearly independant.
Recall that $\Gal(\bL/\bK)=\{\mathfrak{a}_k|
0\le k\le n-1\}$ is generated by $\mathfrak{a}_1$
which satisfies $\mathfrak{a}_1(\theta)=\zeta_n \theta$.
So $$\mathfrak{a}_1(\frac{1}{\zeta_n^{-k} \theta-1})=
\frac{1}{\zeta_n^{-(k-1)} \theta-1}.$$
We conclude that $\Theta$ is a normal basis.
If $i\ne 0$, then
$$\frac{1}{\theta-1}\times \frac{1}{\zeta_n^{-i}\theta-1}=
\frac{1}{\zeta_n^{-i}-1}(\frac{1}{\theta-1}
-\frac{\zeta_n^{-i}}{\zeta_n^{-i}\theta-1}).$$
Since the weight can be defined
using linear combinations in equation $(\ref{eq:30})$,
we conclude that the weight of $\Theta$
is at most $3n-2$.

\subsection{Complexity}\label{section:complexity-multiplicative-group}

We use the same procedure as in section 
\ref{section:complexity-additive-group}. 
The notation is also the same,
except that $$u=\frac{1}{ x-1}, \ \theta_0= \frac{1}{\theta-1},$$
and  for $0 \le k \le p-1$
$$ u_k=u \circ  \tau^k,  \theta_k=
\mathfrak{a}_{-k}(\theta_0), \xi_k=\theta_k^2.$$
We set $\xi_0=\sum_{0\le k\le n-1} \imath_k
\theta_k$ and $\overrightarrow{\imath}=(\imath_k)_{0\le k\le 
n-1}$.  We denote by $R$ the point 
in $\mathbf{G}_m(\bK)$ with coordinate
$x(R)=\zeta_{mn}$. The evaluation map
$$f  \mapsto (f(R+jt))_{0\le k\le n-1}$$
 is a bijection from 
the linear space $\mathcal{L}$ onto $\bK^n$.
Its inverse map is 
$$(a_i)_{0\le i\le n-1}\mapsto \overrightarrow{u}_R^{-1}
\star (a_i)_{0\le i\le n-1},$$
 where
$\overrightarrow{u}_{\!\!R}=( u_0(R+kt))_{0\le k\le n-1}$ and
$\overrightarrow{u}_{\!\!R}^{-1}$ is the inverse of
$\overrightarrow{u}_{\!\!R}$ for the convolution product.

The coordinates in $\Theta$ of the product 
$$\left( \sum_{0\le k\le n-1} \alpha_i\theta_i \right)  \times
\left( \sum_{0\le k\le n-1}\beta_j\theta_j\right)$$
are given by
$$ \overrightarrow{\imath}\star (\overrightarrow{\alpha}\diamond
\overrightarrow{\beta}) +
\overrightarrow{u}_{\!\!R}^{-1}\star \left((\overrightarrow{u}_{\!\!R}\star \overrightarrow{\alpha})
\diamond (\overrightarrow{u}_{\!\!R}\star \overrightarrow{\beta})
- \overrightarrow{w}_{\!\!R}\star (\overrightarrow{\alpha}
\diamond \overrightarrow{\beta})\right),$$
where $\overrightarrow{w}_{\!\!R}
=( u_0^{2}(R+kt))_{0\le k\le n-1}.$This multiplication algorithm 
consists in 5 convolution
products, 2 component-wise products, 1 addition
and 1 subtraction between vectors in $\bK^n.$
Using fast algorithms for the
convolution products, this multiplication algorithm has complexity 
$O(n(\log n)|\log|\log n||)$.

\subsection{Example}

Take $p=61$ and $n=6$, we set $\bK= \bF\!_{61}$.
The point $a \in  \mathbf{G}_m(\bK)$ with coordinate $x(a)=
2 \bmod 61$ generates the group $\mathbf{G}_m(\bK))$, and $48
\bmod 61$ is a primitive $6$-th root of unity in $\bK$.
So $a \bmod I(\mathbf{G}_m(\bK))$ has order $6$ in
$\mathbf{G}_m(\bF\!_{61})/I(\mathbf{G}_m(\bK))$.
Hence $X^6-2\in \bK[X]$ is irreducible over $\bK$. We set
$\bL=\bK[X]/(X^6-2)$ and $\theta=X\bmod X^6-2$.
Thus $\Theta=(\theta_k)_{0\le k\le 5}$ is a normal basis
of $\bL/\bK$.
The weight of $\Theta$ is equal to $15$, according
to  following equations :
$$\theta_0^2 = 53\theta_1+40\theta_2+23\theta_3
+50\theta_4+18\theta_5, \quad \theta_0\times\theta_1= 47(\theta_0-14\theta_1),
\quad \theta_0\times\theta_2=56(\theta_0-13\theta_2),$$

$$\theta_0\times\theta_3= 30(\theta_0+\theta_3),
\quad \theta_0\times\theta_4= 4(\theta_0-47\theta_4),
\quad \theta_0\times\theta_5= 13(\theta_0-48\theta_5).$$

We are going to compute the coordinates 
 in $\Theta$ of the product 
$$\left( \sum_{0\le i \le 5} \alpha_i\theta_i \right)  \times
\left( \sum_{0\le i \le 5}\beta_j\theta_j\right) \text{  where  }
\overrightarrow{\alpha}=(1,3,1,1,2,1) \text{ and } \overrightarrow{\beta}
=(2,1,1,4,2,1).$$

We know that $\overrightarrow{\imath} =(0, 53, 40, 23, 50, 18).$
We compute
$$\overrightarrow{u}_{\!\!a}=(1, 9, 21, 20, 22, 52)$$
$$  \overrightarrow{u}_{\!\!a}^{-1}=
(43, 11, 37, 55, 46, 32), $$
$$\overrightarrow{w}_{\!\!a}=
(1, 20, 14, 34, 57, 20).$$
So
$$
\begin{array}{rl}
 \overrightarrow{\imath}\star (\overrightarrow{\alpha}\diamond
\overrightarrow{\beta}) & =  (43, 29, 46, 36, 12, 32)   \\
  \overrightarrow{u}_{\!\!a}\star \overrightarrow{\alpha}  & = 
 (6, 25, 43, 36, 44, 56)   \\
  \overrightarrow{u}_{\!\!a}\star \overrightarrow{\beta} & = 
   (26, 14, 4, 27, 29, 16) \\
  \overrightarrow{w}_{\!\!a}\star (\overrightarrow{\alpha}
\diamond \overrightarrow{\beta}) & = (26, 14, 4, 27, 29, 16)\\
 (\overrightarrow{u}_{\!\!a}\star \overrightarrow{\alpha}) \diamond 
 (\overrightarrow{u}_{\!\!a}\star \overrightarrow{\beta})  & =
 (52, 45, 50, 57, 56, 42) .
\end{array}$$
Therefore
$$ 
\begin{array}{rl}
 \overrightarrow{u}\!_a^{-1}\star \left( ( \overrightarrow{u}_{\!\!a}\star
 \overrightarrow{\alpha}) \diamond ( \overrightarrow{u}_{\!\!a}\star
 \overrightarrow{\beta}) - \overrightarrow{w}_{\!\!a}\star
 (\overrightarrow{\alpha} \diamond \overrightarrow{\beta})\right)
 & =  (24, 25, 20, 28, 33, 54).
\end{array}$$

Finally we get
$$\left( \sum_{0\le i \le 5} \alpha_i\theta_i \right)  \times
\left( \sum_{0\le j \le 5}\beta_j\theta_j\right) = 6\theta_0 + 54\theta_1
+ 5\theta_2 + 3\theta_3+ 45\theta_4 + 25\theta_5.$$


\section{The Lucas torus case}

This section is devoted to the use of Lucas torus for
constructing normal bases.

\subsection{Basic facts concerning the Lucas Torus}\label{section:Torus-Basic.Facts}

Let $\bK$ be a field with characteristic
different from $2$. Let $\alpha\in\bK$ be a nonsquare element. 
The Lucas torus $\mathbb{T}_\alpha$ over $\bK$ is the affine
plane curve defined by
\begin{equation}\label{eq:11} 
\mathbb{T}_\alpha : x^2-\alpha y^2=1.
\end{equation}
 
This is a commutative algebraic group
with group law $ \oplus_{\mathbb{T}_\alpha}$ defined by
\begin{equation}\label{eq:21} 
(x,y)\oplus_{\mathbb{T}_\alpha}(x',y')=(xx'+\alpha yy',
xy'+x'y).
\end{equation}
Its unit element $O_{\mathbb{T}_\alpha}$ has coordinates equal to $(1, 0)$.

It is easily checked that
the map 
$$
\xymatrix{
 \varphi : \mathbb{T}_\alpha(\bK) & \ar@{->}[r] & \{\xi=(a+b\sqrt{\alpha}) \in \bK(\sqrt{\alpha}) | 
\Norm(\xi)=1\}\\
(x,y) & \ar@{|-{>}}[r] & x+y\sqrt{\alpha}
}$$
is a group isomorphism.

Assume that $\bK = \bF\!_q$ is a finite field with characteristic $p>2$, and
$ \alpha \in \bF\!_q$ is a nonsquare element. We have an exact sequence
\[
\xymatrix{
0 \ar@{->}[r] & \mathbb{T}_\alpha(\bF\!_q) \ar@{->}^{\imath}[r] & \bF\!_{q^2}^*
\ar@{->}^{\Norm}[r] & \bF\!_q^* \ar@{->}[r] & 0
}
\]
Indeed the $\imath$ map is injective because the above
map $\varphi$ is injective. 
Further the $\Norm$
map is surjective [see \cite{Serre-Cours-in-Arithmetic}, proposition 4, page 34].
So $\mathbb{T}_\alpha(\bF\!_q)$ is a cyclic group with order
$(q^2-1)/(q-1)=q+1$.

Actually the map $\varphi :(x,y)\mapsto z:=x+y\sqrt{\alpha}$ is
an isomorphism of $\bK(\sqrt{\alpha})$-varieties,
from the Lucas torus onto
the multiplicative group. The inverse map is given by
$z\mapsto (\frac{z^2+1}{2z},
\frac{z^2-1}{2z\sqrt{\alpha}})$. So $\mathbb{T}_\alpha$ is 
a twist of the multiplicative group.

\subsection{Specializing isogenies of the Lucas torus}\label{section:Basis-Lucas}

Assume that $n$ is a non-trivial divisor of $q+1$.
We denote by $m$ the cofactor of $n$, that is $q+1=nm$.
The $n$-torsion $\mathbb{T}_\alpha[n]$ is a cyclic subgroup
of $\mathbb{T}_\alpha(\bF\!_q)$.
Let $\Phi : \overline{\bF}\!_q \rightarrow 
\overline{\bF}\!_q$ be the $q$-Frobenius automorphism.
Let $I : \mathbb{T}_ \alpha
\rightarrow \mathbb{T}_ \alpha$ be the quotient isogeny
of $\mathbb{T}_ \alpha$ by $\mathbb{T}_\alpha[n]$.
This is the multiplication by $n$ isogeny. 
We have
$$\begin{array}{l}
{\varphi( I(x,y))=\varphi(x,y)^n=(x+\sqrt{\alpha}y)^n}\\
{= \sum_{0\le 2k\le 2n} \left(\begin{array}{l}
n\\
2k \end{array} \right) x^{n-2k} y^{2k}\alpha^k +
\sqrt{\alpha} \sum_{0\le 2k+1\le 2n} \left(\begin{array}{l}
n\\
2k+1 \end{array} \right) x^{n-2k-1} y^{2k+1}\alpha^k}.
\end{array}$$
The quotient group $\mathbb{T}_\alpha(\bF\!_q)/I(\mathbb{T}_\alpha(\bF\!_q))$
is cyclic of order $n$. Let $a$ be a generator of $\mathbb{T}_\alpha
(\bF\!_q)$. Then $a\bmod I(\mathbb{T}_\alpha(\bF\!_q))$ generates
$\mathbb{T}_\alpha(\bF\!_q)/I(\mathbb{T}_\alpha(\bF\!_q))$.
The subfield $\bL = \bF\!_q(I^{-1}(a))$ of $\overline{\bF}\!_q$
is a cyclic extension of $\bF\!_q$ with degree $n$.
Indeed $\bL$ is separable because it is an algebraic extension
of the perfect field $\bF\!_q$.
The fiber $I^{-1}(a)$ is defined over $\bF\!_q$ because
our Lucas torus is defined over
$\bF\!_q$ and $a$ is a $\bF\!_q$-rational point. Thus
any $\bF\!_q$-automorphism
of $\overline{\bF}\!_q$ maps $I^{-1}(a)$ into itself.
So the $\bF\!_q$-automorphisms
of $\overline{\bF}\!_q$ map $\bL$ into
itself, that is $\bL$ is normal over $\bF\!_q$.
Let $b$ be a point in $I^{-1}(a)$.
For any $\sigma \in \Gal(\bL/\bF\!_q)$ we have
\begin{equation}\label{eq:12} 
^\sigma b\ominus_{\mathbb{T}_\alpha}b \in \Ker(I),
\end{equation}
that is there exists $t_\sigma \in \mathbb{T}_\alpha[n]$
such that
\begin{equation}\label{eq:13} 
^\sigma b=b\oplus_{\mathbb{T}_\alpha}t_\sigma.
\end{equation}
If $\sigma, \sigma'$ lie in $\Gal(\bL/\bF\!_q)$ then
$$^{\sigma \circ \sigma'}\! b= \ ^{\sigma}\!(b\oplus_{\mathbb{T}_\alpha}t_{\sigma'})=
b\oplus_{\mathbb{T}_\alpha}t_{\sigma'}
\oplus_{\mathbb{T}_\alpha}t_{\sigma}.$$
So $^{\sigma \circ \sigma'} b= \  ^{\sigma' \circ \sigma}b$.
Hence $\Gal(\bL/\bF\!_q)$ is an abelian group
with exponent dividing $n$.
Let $\langle a\bmod I(\mathbb{T}_\alpha(\bF\!_q))\rangle$ be
the subgroup of $\mathbb{T}_\alpha(\bF\!_q)/I(\mathbb{T}_\alpha(\bF\!_q))$
generated by the class of $a$.
The map
$$\xymatrix{
\kappa : \left( \langle a\bmod I(\mathbb{T}_\alpha(\bF\!_q))\rangle/I(\mathbb{T}_\alpha(\bF\!_q)) \right) 
\times \Gal(\bL/\bF_q) \ar@{->}[r] & \mathbb{T}_\alpha[n]\\
(k\overline{a},\sigma) \ar@{|->}[r] & (kb)^\sigma -(kb),
}$$
is a non-degenerate pairing. It induces a group isomorphism
$$\kappa_\bL : \Gal(\bL/\bF\!_q)
\rightarrow  \Hom( \langle a\bmod I(\mathbb{T}_\alpha(\bF\!_q))\rangle/I(\mathbb{T}_\alpha(\bF\!_q)) ,
\mathbb{T}_\alpha[n]).$$
The group
 $\Hom( \langle a\bmod I(\mathbb{T}_\alpha(\bF_q))\rangle/I(\mathbb{T}_\alpha(\bF\!_q)), 
\mathbb{T}_\alpha[n])$ has order $n$.
So $\#\Gal(\bL/ \bF\!_q)= n$, that is $\bL$ is a degree $n$ cyclic
extension of $\bF\!_q$ with Galois group generated by 
the $q$-Frobenius automorphism $\Phi$.
We set 
$t=b^\Phi \ominus_{\mathbb{T}_\alpha}b$.
This a generator of the $n$-torsion
$\mathbb{T}_\alpha[n]$.
So the fiber of $I$ above $a$ is given by 
\begin{equation}\label{eq:15}
I^{-1}(a)=[b]+[b\oplus_{\mathbb{T}_\alpha}t]+\cdots+[b
\oplus_{\mathbb{T}_\alpha}(n-1)t)].
\end{equation}

The projective closure of $\mathbb{T}_\alpha$
is the locus
\begin{equation}\label{eq:20}
\overline{\mathbb{T}}_\alpha : X^2-\alpha Y^2-Z^2=0
\end{equation}
in the projective plane $\cP^2$.
This is a genus $0$ projective curve with two points
$\infty_1=(\sqrt{\alpha}:1:0)$ and $\infty_2=(-\sqrt{\alpha}:1:0)$ on the
line at infinity. The translation
$\tau : P \mapsto P\ominus_{\mathbb{T}_\alpha}t$ is an automorphism of the
Lucas torus which extends
to its projective closure $\overline{\mathbb{T}}_\alpha$
by setting 
$\tau(\infty_1)=\infty_1$, and  $\tau(\infty_2)=\infty_2$.

Denote by $\Div^0(\mathbb{T}_\alpha)$ the subgroup
of $\Div^0(\overline{\mathbb{T}}_\alpha)$ made of divisors
$D$ with support contained in
$\mathbb{T}_\alpha(\overline{\bF}\!_q)$.
A divisor $D$ in $\Div^0(\mathbb{T}_\alpha)$ is said to be principal if
there exists a non-zero function $f$
in $\overline{\bF}\!_q(\overline{\mathbb{T}}_\alpha)$
such that $f(\infty_1)=f(\infty_2)$ and $D=\div(f)$.

\noindent \textbf{Claim.}
If $D=\div(f)$ is a principal divisor in $\Div^0(\mathbb{T}_\alpha)$,
then  $\sum_{P \in \overline{\mathbb{T}}_\alpha }
 [\ord_P(f)]P=O_{\mathbb{T}_\alpha}$ and $f$ has at least two poles.

\noindent Indeed,
using arguments similar to that used in [\cite{Koblitz}, section 4],
one shows that a point $P $ in $\mathbb{T}_\alpha$ is a zero (resp. a pole)
of $f$ if and only if $-P$ is a zero
(resp. a pole) of $f$. More precisely, for any $P$
in $\mathbb{T}_\alpha$, we have $\ord_P(f)=\ord_{-P}(f)$.
So $$\sum_{P \in \overline{\mathbb{T}}_\alpha }
 [\ord_P(f)]P=O_{\mathbb{T}_\alpha}.$$
Since $\div(f)$ has degree zero, we are done.

Set 
$$x = X/Z, \  y = Y/Z \ \text{ and }  \ v = \frac{x-1}{y}.$$
Using 
equation $(\ref{eq:20})$, we obtain
$$y=\frac{(x+1)v}{\alpha}, \  x= \frac{(x+1)v^{2}}{\alpha}+1.$$

Besides,
$$\div(v)= [O_{\mathbb{T}_\alpha}]-  [\widetilde{O}_{\mathbb{T}_\alpha}],$$
where $\widetilde{O}_{\mathbb{T}_\alpha}$ is the
 point in $\mathbb{T}_\alpha$ with coordinate $(-1,0)$.
Let $u_{O,t}$ be the function in $\bF_q(\overline{\mathbb{T}}_\alpha)$
defined by
$$u_{O,t}=1+\frac{1}{y-\frac{y(t)}{x(t)-1}(x-1)}.$$
Its polar divisor 
is equal to $ - [O_{\mathbb{T}_\alpha}]-[ t ],$ and we have
$u_{O,t}(\infty_1)=u_{O,t}(\infty_2)=1$.\\
For $0\le k\le n-1$, we set 
\begin{equation}\label{eq:16}
y_{k}=y\circ \tau^{k}, \ x_{k}=x\circ \tau^{k}, \ v_k=v \circ \tau^{k}
 \text{ and } u_{kt,(k+1)t}:=u_{O,t} \circ \tau^{k}.
\end{equation}
The function $v_k$ is a uniformizer at $kt$. Besides,
the function $u_{kt,(k+1)t}$ has two simple poles ($kt$ and $(k+1)t$),
and $\div(u_{kt,(k+1)t})$ is a principal divisor on $\mathbb{T}_\alpha$.
Since
$$u_{kt,(k+1)t}=1+\frac{1}{y_{k}} \times \frac{1}{1-\frac{y(t)}{x(t)-1} \times v_{k}},$$
the Taylor expansion of $u_{kt,(k+1)t}$ at $kt$  is
\begin{equation}\label{eq:17}
u_{kt,(k+1)t}=1 + \frac{\alpha }{2}(\frac{1 }{v_k}+\frac{ y(t)}{(x(t)-1)}+\frac{y^2(t)}{(x(t)-1)^2}v_k+ O(v_k^2)).
\end{equation}
Since $$u_{kt,(k+1)t}=
1+\frac{1}{y(P-(k+1)t+t)-\frac{y(t)}{x(t)-1}(x(P-(k+1)t+t)-1+x_{k+1}-x_{k+1})}$$
$$=1-\frac{1}{y_{k+1}} \times \frac{1}{1+\frac{y(t)}{x(t)-1} \times v_{k+1}},$$
the Taylor expansion of $u_k$ at $(k+1)t$ is 
$$u_{kt,(k+1)t}=1 - \frac{\alpha }{2}(\frac{1}{v_{k+1}} - \frac{ y(t)}{ (x(t)-1)}+
\frac{ y^2(t)}{(x(t)-1)^2}v_{k+1}- O(v_{k+1}^2)).$$

The divisor
$$D=[0_{\mathbb{T}_\alpha}]+[t]+[2t]+
\cdots+[(n-1)t].$$ 
induces a linear space
$$\mathcal{L}=\mathbf{H}^0(\overline{\mathbb{T}}_\alpha,
\mathcal{O}_{\overline{\mathbb{T}}_\alpha}
(D))$$ with dimension $n+1$.
Denote by $\mathcal{L}_1$ the subspace of $\mathcal{L}$
made of functions $f\in \mathcal{L}$ such that $f(\infty_1)=f(\infty_2)$.
This  is a $n$-dimensional $\overline{\bF}\!_q$-vector space, as shown by the
following lemma.

\begin{lemma}\label{lemma1}
With the above notation:
\begin{enumerate}
\item The sum $\sum_{0\le k\le n-1} u_{kt,(k+1)t}$ is a constant
$\mathfrak{c} \in \bF\!_q$, and there exists two scalars 
$\mathfrak{a}\ne 0$ and $\mathfrak{b}$ in $\bF\!_q$
such that $ \mathfrak{a}\mathfrak{c} +n\mathfrak{b}=1.$

\item We set $u_k=\mathfrak{a} u_{kt,(k+1)t} +\mathfrak{b}$.
Then $(u_0,u_1,\ldots,
u_{n-1})$ is a basis of $\mathcal{L}_1$.
\item The system $\Theta=(\theta_0,\theta_1,\ldots,
\theta_{n-1})$ defined by $\theta_k=u_k(b)$
is a normal basis of $\bL$ over $\bF\!_q$.
\end{enumerate}
\end{lemma}

\begin{myproof}
\begin{enumerate}
\item Since $\div(\sum_{0\le k\le n-1} u_{kt,(k+1)t})$
 is a principal divisor on $\mathbb{T}_\alpha$,
the proof uses the arguments of [\cite{Couveignes-Lercier2}, proof of lemma 4].
First we denote by $T$ the subgroup generated by $t$.
The sum $\sum_{0\le k\le n-1} u_{kt,(k+1)t}$ is invariant by translations
in $T$. So it can be seen as a function on the quotient of
$\mathbb{T}_\alpha$ by $T$, with no more than one pole and taking
the same value at two distinct points.
Hence it is a constant which we denote by $\mathfrak{c}$.
If the characteristic $p$ of $\bF\!_q$ divides $n$,
it is easily checked that $\mathfrak{c}\ne 0$ using the
function $\sum_{0\le k\le n-1} k u_{kt,(k+1)t}$.
So either $n$ is prime to $p$ or $\mathfrak{c}\ne 0$.
In any case, there exists two scalars 
$\mathfrak{a}\ne 0$ and $\mathfrak{b}$ in $\bF\!_q$
such that $ \mathfrak{a}\mathfrak{c} +n\mathfrak{b}=1.$

\item and 3. Let $\lambda_0,\lambda_1,\ldots,\lambda_{n-1}$ be scalars in
$\bF\!_q$ such that
$$\sum_{k\in \ZZ/n\ZZ}\lambda_k \theta_k=0.$$
Then the function
$F=\sum_{k\in \ZZ/n\ZZ}\lambda_k u_k \in \bF\!_q[X]$ 
cancels at $b$ and also at all its $n$ conjugates over $\bF\!_q$.
Since $\div(\sum_{k\in \ZZ/n\ZZ}\lambda_k u_k)$ is a principal
divisor on $\mathbb{T}_\alpha$ and $b$ does not lie
in the $n$-torsion, we deduce that
$\sum_{k\in \ZZ/n\ZZ}\lambda_k u_k=0$.
Examination of
poles shows that all $\lambda$'s are equal.
So
$$\lambda_0(\sum_{k\in \ZZ/n\ZZ} u_k)=0.$$
But $\sum_{k\in \ZZ/n\ZZ} u_k=1$. So all
$\lambda$'s are equal to zero.

Recall that $\Gal(\bL/\bF\!_q)$ is generated by $\Phi$
which satisfies $b^\Phi=b+t$. So
$\Phi(\theta_k)=\theta_{k-1}$.
\end{enumerate}

\end{myproof}

If $k,l \in \ZZ/n\ZZ$ and $k\ne l, l+1,l-1 \bmod n$, then
$$u_ku_l \in \mathcal{L}_1.$$
Further
$$u_{k-1}u_k+\frac{\alpha^2}{4v^2_{k}} \in \mathcal{L}_1,$$
 and 
$$ u_k^2-\frac{\alpha^2}{4v^2_{k}}-
 \frac{\alpha^2}{4v^2_{k+1}} \in \mathcal{L}_1.$$

If $\sum_k a_k u_k$ and $\sum_k b_k  u_k$ are two functions
in $\mathcal{L}_1$, then
$$
(\sum_k a_k u_k)(\sum_k b_k 
u_k) = \frac{\alpha^2}{4}\sum_k a_k b_k
\left(\frac{1}{v^2_k}
+\frac{1}{v^2_{k+1}} \right) -
\frac{\alpha^2}{4}\sum_k 
\frac{a_{k-1}b_k}{v^2_k} 
-\frac{\alpha^2}{4}\sum_k 
\frac{b_{k-1}a_k}{v^2_k} \bmod 
\mathcal{L}_1,$$
that is
\begin{equation}\label{eq:18}
(\sum_k a_k u_k)(\sum_k b_k u_k)
=\frac{\alpha^2}{4}\sum_k
\frac{(a_k-a_{k-1})
(b_k - b_{k-1})}{v^2_k} 
 \bmod 
\mathcal{L}_1.
\end{equation}

\subsection{Complexity}\label{section:complexity-Lucas torus}

We use the same procedure as in section 
\ref{section:complexity-additive-group} . The notation is also
the same, except that
$$u= \mathfrak{a}(1+\frac{1}{y-\frac{y(t)}{x(t)-1}(x-1)})+ \mathfrak{b}, \
v=\frac{x-1}{y}, $$
and for $0 \le k \le n-1$ we set
$$u_k=u\circ \tau^k, \ \theta_k=u_k(b), \ v_k=v\circ \tau^k ,  \  \xi_k=
\frac{1}{v_k^2(b)}.$$
We also set $\xi_0=\sum_{0\le k\le n-1} \imath_k
\theta_k$ and $\overrightarrow{\imath}=(\imath_k)_{0\le k\le 
n-1}$. A linear combination of the $\xi$'s can be reduced into a linear combination
of the $\theta$'s
$$\sum_{0\le i \le n-1} x_i \xi_i =
\sum_{0\le k \le n-1} (\overrightarrow{\imath}\star_k \overrightarrow{x}) \theta_k,$$
where $\overrightarrow{x}=(x_0,x_1,\ldots, x_{n-1})$.
The evaluation map $f\mapsto (f(a+jt)_{0\le j \le n-1}$
is a bijection from the linear space $\mathcal{L}_1$ onto $\bF_q^n$
because $[n](a)\ne O_{\mathbb{T}_\alpha}$.
Its inverse map is 
$$\overrightarrow{x} \mapsto \overrightarrow{u}_{\!\!a}^{-1}\star\overrightarrow{x},$$
where $\overrightarrow{u}_{\!\!a}^{-1}$ is the inverse 
of $\overrightarrow{u}_{\!\!a}=(u_0(a+jt))_{0\le j\le n-1 }$ for the convolution product.

So the coordinates in $\Theta$ of a product 
$$\left( \sum_{0\le i \le n-1} x_i\theta_i \right)  \times
\left( \sum_{0\le j \le n-1} y_j\theta_j\right)$$
are given by
$$
\begin{array}{ll}
 (\frac{\alpha^2}{4} \overrightarrow{\imath})\star \left( (\overrightarrow{x} -\sigma(\overrightarrow{x}))\diamond
(\overrightarrow{y}-\sigma(\overrightarrow{y}))\right) + & \\
 & \\
 \overrightarrow{u}_{\!\!a}^{-1}\star \left[(\overrightarrow{u}_{\!\!a}\star \overrightarrow{x})
\diamond (\overrightarrow{u}_{\!\!a}\star  \overrightarrow{y})
- (\frac{\alpha^2}{4} \overrightarrow{w}_{\!\!a})\star \left( \left(\overrightarrow{x} -
\sigma(\overrightarrow{x})\right)\diamond
(\overrightarrow{y}-\sigma(\overrightarrow{y}))\right)\right],
& \\
\end{array}
$$
where $\overrightarrow{y}=(y_0,y_1,\ldots, y_{n-1})$ and $\overrightarrow{w}_{\!\!a}=(\frac{1}{v_0^2(a+jt)})_{0\le j\le n-1 }$.
\noindent That consists in 5 convolution
products, 2 component-wise products, 1 additions
and 3 subtraction between vectors in $\bF_q^n$.
So this multiplication algorithm has complexity 
$O(n(\log n)|\log|\log n||)$.

\subsection{Example}

Take $p=q=7$, $n=4$, and $\alpha=3$. We set $\bK= \bF_{7}$. 
The point $a=(5,1)$ generates the group $\mathbb{T}_3(\bF_{7})$.
 The point $t=(0,3)$ generates the $4$-torsion subgroup.
 The multiplication by $4$ isogeny $I : \mathbb{T}_3 
 \rightarrow  \mathbb{T}_3$ is given in terms of coordinates by
 $$I(x,y)= (x^4 + 4y^2x^2 +2y^4, 4yx^3 + 5y^3x).$$
We focus on the intersection points of the two planes curves 
$$\{ x^4 + 4y^2x^2 +2y^4 -5=0\} \ \text{ and } \
\{ 4yx^3 + 5y^3x -1 =0\}$$
Since their defining polynomials are relatively prime,
they have finitely many intersection points.
The polynomial $P(x)=x^4+x^2+3$ is an irreducible
factor of the resultant of the defining polynomials.
 We set $\bL=\bK[x]/(x^4 + x^2 + 3)$,
 $\theta=x\bmod x^4 + x^2 + 3$, and 
 $b=(\theta,6\theta^3+5\theta)$. 

We find 
$$(u_{kt,(k+1)t})_{k\in \ZZ/4\ZZ}=(\frac{1}{y+3(x-1)},\frac{-1}{3(x+9y+1)},
\frac{1}{3(9y+27x-1)}, \frac{-1}{3(27x+243y+1)}),$$
so that $\mathfrak{c}=0, \mathfrak{a}=1, \mathfrak{b}=2$ and

$$\theta_0 = 2\theta^3+6\theta^2+\theta+1, \ 
\theta_1=\theta^3+\theta^2+\theta+2,$$
$$\theta_2=5\theta^3+6\theta^2+6\theta+1, \ 
\theta_3=6\theta^3+\theta^2+6\theta+2.$$

We are going to compute the coordinates 
 in $\Theta=(\theta_k)_{0\le k\le 3}$ of the product 
$$\left( \sum_{0\le i \le 3} x_i\theta_i \right)  \times
\left( \sum_{0\le j \le 3} y_j\theta_j\right) \text{  with  }
\overrightarrow{x}=(1,3,1,1)
et \text{ and } \overrightarrow{y}
=(2,1,1,4).$$

We know that 
$$\overrightarrow{\imath} =(3, 0, 0, 3).$$
We compute
$$\overrightarrow{u}_{\!\!a}=(1, 4, 4, 0)$$
$$  \overrightarrow{u}_{\!\!a}^{-1}=
(0, 2, 6,  3), $$
$$\overrightarrow{w}_{\!\!a}=
(1, 2, 2, 0).$$
So
$$
\begin{array}{rl}
 (\frac{\alpha^2}{4} \overrightarrow{\imath})\star \left( (\overrightarrow{x} -\sigma(\overrightarrow{x}))\diamond
(\overrightarrow{y}-\sigma(\overrightarrow{y}))\right) & =  (4, 4, 0, 0)   \\
  \overrightarrow{u}_{\!\!a}\star \overrightarrow{\alpha}  & = 
 (2, 4, 3, 3)   \\
  \overrightarrow{u}_{\!\!a}\star \overrightarrow{\beta} & = 
   (1, 4, 6, 5) \\
  (\frac{\alpha^2}{4} \overrightarrow{w}_{\!\!a})\star \left( (\overrightarrow{x} -\sigma(\overrightarrow{x}))\diamond
(\overrightarrow{y}-\sigma(\overrightarrow{y}))\right) & = (0, 6, 5, 5)\\
 (\overrightarrow{u}_{\!\!a}\star \overrightarrow{\alpha}) \diamond 
 (\overrightarrow{u}_{\!\!a}\star \overrightarrow{\beta})  & =
 (2, 2, 4, 1) .
\end{array}$$
Therefore
$$ 
\begin{array}{rl}
 \overrightarrow{u}_{\!\!a}^{-1}\star \left( ( \overrightarrow{u}_{\!\!a}\star
 \overrightarrow{\alpha}) \diamond ( \overrightarrow{u}_{\!\!a}\star
 \overrightarrow{\beta}) - (\frac{\alpha^2}{4} \overrightarrow{w}_{\!\!a})\star \left( (\overrightarrow{x} -\sigma(\overrightarrow{x}))\diamond
(\overrightarrow{y}-\sigma(\overrightarrow{y}))\right) \right)
 & =  (2, 3, 6, 3).
\end{array}$$

Finally we get
$$\left( \sum_{0\le i \le 3} \alpha_i\theta_i \right)  \times
\left( \sum_{0\le j \le 3}\beta_j\theta_j\right) = 6\theta_0 
+ 6\theta_2 + 3\theta_3.$$


\vspace{2cm}

\begin{thebibliography}{2}


\bibitem{Ash-Blake-Vanstone}
David W. Ash, Ian F. Blake and Scott A. Vanstone
\newblock  Low complexity normal bases.
\newblock {\em  Discrete Appl. Math.},  25 (1989), no. 3, pp. 191-210.


\bibitem{Cantor-Kaltofen}
David G. Cantor and Erich Kaltofen
\newblock On fast multiplication of polynomials over arbitrary algebras
\newblock {\em Acta Inform. } 28 (1991), no. 7, pp. 693-701.



\bibitem{Christopoulou-Garefalakis-Panario-Tomson1}
M. Christopoulou, T. Garefalakis, D. Panario and D.Thomson
\newblock The trace of an optimal normal element and low complexity normal bases.
\newblock {\em Des. Codes Cryptogr.},  49 (2008), no. 1, pp. 199-215.




\bibitem{Couveignes-Lercier2}
Jean-Marc Couveignes and Reynald Lercier.
\newblock Elliptic periods for finite fields.
\newblock {\em	Finite Fields Appl.},  15 (2009), no. 1, pp. 1-22.


\bibitem{Gao-Phd}
Shuhong Gao.
\newblock Normal bases over finite fields.
\newblock {\em Thesis (Ph.D.)--University of Waterloo (Canada)},  1993

\bibitem{Gao-Lenstra}
Shuhong Gao and Jr. Hendrik W. Lenstra
\newblock Optimal normal bases.
\newblock {\em Des. Codes Cryptogr.}, 2 (1992), no. 4,pp. 315-323

\bibitem{Gao-Gathen-Panario-Shoup}
S. Gao, J. von zur Gathen, D. Panario and V.Shoup
\newblock Algorithms for exponentiation in finite fields.
\newblock {\em J. Symbolic Comput.}, (2000), pp. 879-889


\bibitem{Koblitz}
Neal Koblitz,
\newblock Algebraic aspects of cryptography.
\newblock {\em Algorithms and Computation in Mathematics, Springer-Verlag, Berlin},
1998.

\bibitem{Liao-You},
Qunying Liao and Lin You.
\newblock Low complexity of a class of normal bases over finite fields.
\newblock {\em Finite Fields and Their Applications}, 2011, pp. 1-14

\bibitem{milne2006}
James Stuart Milne
\newblock  Elliptic Curves.
\newblock {\em BookSurge Publishers}, 2006

\bibitem{Mullin-al}
R.C. Mullin, I.M.  Onyszchuk, S.A. Vanstone, and R.C. Wilson
\newblock  Optimal normal bases $GF(p^n)$.
\newblock {\em  Discrete Applied Math.}, 22 (1988/1989), pp. 149-161



\bibitem{schonhage-Strassen}
A. Sch\"onhage  and  V. Strassen
\newblock  Schnelle {M}ultiplikation grosser {Z}ahlen.
\newblock {\em Computing (Arch. Elektron. Rechnen)}, 7 (1971), pp. 281-292


\bibitem{schonhage}
A. Sch\"onhage 
\newblock Schnelle {M}ultiplikation von {P}olynomen \"uber {K}\"orpern der
{C}harakteristik 2.
\newblock {\em Acta Informat.}, 7 (1976/77), no. 4, pp. 395-398


\bibitem{Serre-Cours-in-Arithmetic}
Jean-Pierre Serre.
\newblock A course in arithmetic.
\newblock {\em Graduate Texts in Mathematics, No. 7, Springer-Verlag, New York-Heidelberg}, 1973

\bibitem{Gathen-Gerhard}
Joakim von zur Gathen and J\"urgen Gerhard.
\newblock Modern computer algebra.
\newblock {\em Cambridge University Press, Cambridge}, 2013


\bibitem{Wan-Zhou}
Zhe-Xian Wan and Kai Zhou.
\newblock On the complexity of the dual basis of a type {I} optimal normal basis.
\newblock {\em Finite Fields Appl.}, 13 (2007),  no. 2, pp. 411-417

\bibitem{Wassermann1}
Alfred Wassermann.
\newblock Konstruktion von {N}ormalbasen.
\newblock {\em  Bayreuther Mathematische Schriften}, 31 (1990), pp. 155-164



\end{thebibliography}

\end{document}